\documentclass{amsart}

\usepackage{fancyhdr}
\usepackage{lastpage}
\usepackage{stmaryrd,yhmath}

\newcommand{\bdis}{\begin{displaymath}}
\newcommand{\edis}{\end{displaymath}}
\newcommand{\be}{\begin{equation}}
\newcommand{\ee}{\end{equation}}
\newcommand{\mbb}{\mathbb}
\newcommand{\mcal}{\mathcal}

\newcommand{\vp}{\varphi}

\newcommand{\mG}{\mathring{G}}
\newcommand{\mT}{\mathring{T}}
\newcommand{\mg}{\mathring{g}}

\newcommand{\zf}{\zeta\left(\frac{1}{2}+it\right)}

\theoremstyle{definition}

\theoremstyle{remark}
\newtheorem{remark}[]{Remark}

\newtheorem*{mydef1}{{\bf Theorem}}

\newtheorem*{mydef4}{{\bf Corollary}}

\newtheorem*{mydef5}{{\bf Lemma}}

\numberwithin{equation}{section}

\begin{document}

\title{Jacob's ladders and the nonlocal interaction of the function $Z^2(t)$ with the function $\tilde{Z}^2(t)$ on the distance
$\sim (1-c)\pi(t)$ for the collections of disconnected sets}

\author{Jan Moser}

\address{Department of Mathematical Analysis and Numerical Mathematics, Comenius University, Mlynska Dolina M105, 842 48 Bratislava, SLOVAKIA}

\email{jan.mozer@fmph.uniba.sk}


\begin{abstract}
It is shown in this paper that there is a fine correlation of the fourth order between the functions $Z^2[\vp_1(t)]$ and $\tilde{Z}^2(t)$, respectively.
This correlation is with respect to two collections of disconnected sets. Corresponding new asymptotic formulae cannot be obtained within known theories
of Balasubramanian, Heath-Brown and Ivic.
\end{abstract}

\maketitle

\section{The result}

In this paper we obtain some new properties of the signal
\bdis
Z(t)=e^{i\vartheta(t)}\zf
\edis
generated by the Riemann zeta-function, where
\be \label{1.1}
\vartheta(t)=-\frac t2\ln\pi+\text{Im}\ln\Gamma\left(\frac 14 +i\frac t2\right)=\frac t2\ln\frac{t}{2\pi}-\frac t2-\frac{\pi}{8}+
\mcal{O}\left(\frac 1t\right) .
\ee
Let (see \cite{3})
\be \label{1.2}
\begin{split}
& G_3(x)=G_3(x;T,U)= \\
& =\bigcup_{T\leq g_{2\nu}\leq T+U}\{ t:\ g_{2\nu}(-x)\leq t\leq g_{2\nu}(x)\}, \ 0<x\leq \frac{\pi}{2} , \\
& G_4(y)=G_4(y;T,U)= \\
& =\bigcup_{T\leq g_{2\nu+1}\leq T+U}\{ t:\ g_{2\nu+1}(-y)\leq t\leq g_{2\nu+1}(y)\}, \ 0<y\leq \frac{\pi}{2} ,
\end{split}
\ee
\be \label{1.3}
U=T^{5/12}\ln^3T ,
\ee
and the collection of sequences $\{ g_\nu(\tau)\}$, $\tau\in [-\pi,\pi]$, $\nu=1,2,\dots $ is defined by the equation (see \cite{2}, \cite{3}, (6))
\bdis
\vartheta_1[g_\nu(\tau)]=\frac{\pi}{2}\nu+\frac{\tau}{2};\ g_\nu(0)=g_\nu ,
\edis
where (comp. (\ref{1.1})
\bdis
\vartheta_1(t)=\frac t2\ln\frac{t}{2\pi}-\frac t2-\frac{\pi}{8}.
\edis
Let
\be \label{1.4}
G_3(x)=\vp_1[\mG_3(x)],\ G_4(y)=\vp_1[\mG_4(y)] ,
\ee
where $y=\vp_1(T),\ T\geq T_0[\vp_1]$ is the Jacob's ladder. The following theorem holds true.

\begin{mydef1}
\be \label{1.5}
\begin{split}
& \int_{\mG_3(x)}Z^2[\vp_1(t)]\tilde{Z}^2(t){\rm d}t= \\
& \frac x\pi U\ln\frac{T}{2\pi}+\frac{2x}{\pi}\left( c + \frac{\sin x}{x}\right)U+\mcal{O}(xT^{5/12}\ln^2T), \\
& \int_{\mG_4(y)}Z^2[\vp_1(t)]\tilde{Z}^2(t){\rm d}t= \\
& \frac y\pi U\ln\frac{T}{2\pi}+\frac{2y}{\pi}\left( c - \frac{\sin y}{y}\right)U+\mcal{O}(yT^{5/12}\ln^2T),
\end{split}
\ee
where
\be \label{1.6}
t-\vp_1(t)\sim (1-c)\pi(t),\ t\to\infty ,
\ee
and $c$ is the Euler's constant and $\pi(t)$ is the prime-counting function.
\end{mydef1}

Let (comp. (\ref{1.4})) $T=\vp_1(\mT),\ T+U=\vp_1(\widering{T+U})$. Similarly to \cite{14}, (1.8) we have
\be \label{1.7}
\rho\{ [T,T+U];[\mT,\widering{T+U}]\}\sim (1-c)\pi(T) ,
\ee
where $\rho$ stands for the distance of the corresponding segments.

\begin{remark}
Some nonlocal interaction of the functions $Z^2[\vp_1(t)],\ \tilde{Z}^2(t)$ is expressed by eq. (\ref{1.5}). This interaction is connected with two
collections of disconnected sets unboundedly receding each from other (see (\ref{1.6}), (\ref{1.7});\ $\rho\to\infty$ as $T\to\infty$) - like mutually
receding galaxies (the Hubble law).
\end{remark}

\begin{remark}
The asymptotic formulae (\ref{1.5}) (comp. (1.3)) cannot be received by methods of Balasubramanian, Heath-Brown and Ivic (comp. \cite{1}).
\end{remark}

This paper is a continuation of the series of works \cite{4}-\cite{18}.

\section{On big asymmetry in the behaviour of the function $Z^2[\vp_1(t)]\tilde{Z}^2(t)$ relatively to the sets $\mG_3(x)$ and $\mG_4(x)$}

We obtain from (\ref{1.5})

\begin{mydef4}
\be \label{2.1}
\begin{split}
& \int_{\mG_3(x)}Z^2[\vp_1(t)]\tilde{Z}^2(t){\rm d}t-\int_{\mG_4(x)}Z^2[\vp_1(t)]\tilde{Z}^2(t){\rm d}t= \\
& \frac{4}{\pi}U\sin x+\mcal{O}(xT^{5/12}\ln^2T),\ x\in (0,\pi/2] ,
\end{split}
\ee
especially, in the case $x=\pi/2$, we have
\be \label{2.2}
\int_{\mG_3(\pi/2)}Z^2[\vp_1(t)]\tilde{Z}^2(t){\rm d}t-\int_{\mG_4(\pi/2)}Z^2[\vp_1(t)]\tilde{Z}^2(t){\rm d}t\sim \frac 4\pi U,
\ee
where $[\mT,\widering{T+U}]\subset \mG_3(\pi/2)\cup\mG_4(\pi/2)$.
\end{mydef4}

\begin{remark}
The formulae (\ref{2.1}), (\ref{2.2}) represent the big difference of the areas (measures) of the figures which correspond to the functions
\bdis
Z^2[\vp_1(t)]\tilde{Z}^2(t),\ t\in\mG_3(x) ; \qquad Z^2[\vp_1(t)]\tilde{Z}^2(t),\ t\in\mG_4(x) .
\edis
The reason for this big asymmetry given by (\ref{2.2}) is probably the fact that the zeroes of
\bdis
\zf,\ t\in\mG_3\left(\frac{\pi}{2}\right)\bigcup \mG_4\left(\frac{\pi}{2}\right)
\edis
lie preferably in the set $\mG_4(\pi/2)$.
\end{remark}

\section{Proof of the Theorem}

\subsection{}

The following lemma holds true (see \cite{6}, (2.5); \cite{7}, (3.3); \cite{14}, (4.1))

\begin{mydef5}
For every integrable function (in the Lebesgue sense) $f(x),\ x\in [\vp_1(T),\vp_1(T+U)]$ the following is true
\be \label{3.1}
\int_T^{T+U}f[\vp_1(t)]\tilde{Z}^2(t){\rm d}t=\int_{\vp_1(T)}^{\vp_1(T+U)}f(x){\rm d}x,\ U\in \left(\left. 0,\frac{T}{\ln T}\right]\right. ,
\ee
where $t-\vp_1(t)\sim (1-c)\pi(t)$, $c$ is the Euler's constant, $\pi(t)$ is the prime-counting function and
\bdis
\begin{split}
& \tilde{Z}^2(t)=\frac{{\rm d}\vp_1(t)}{{\rm d}t},\ \vp_1(t)=\frac 12\vp(t), \\
& \tilde{Z}^2(t)=\frac{Z^2(t)}{2\Phi^\prime_\vp[\vp(t)]}=\frac{Z^2(t)}{\left\{ 1+\mcal{O}\left(\frac{\ln\ln t}{\ln t}\right)\right\}\ln t} ,
\end{split}
\edis
(see \cite{7}, (1.1), (3.1), (3.2)).
\end{mydef5}

\begin{remark}
The formula (\ref{3.1}) remains true also in the case when the integral on the right-hand side of (\ref{3.1}) is only relatively convergent integral
of the second kind (in the Riemann sense).
\end{remark}

In the case (comp. (\ref{1.4}) $T=\vp_1(\mT),\ T+U=\vp_1(\widering{T+U})$) we obtain from (\ref{3.1})
\be \label{3.2}
\int_{\mT}^{\widering{T+U}}f[\vp_1(t)]\tilde{Z}^2(t){\rm d}t=\int_T^{T+U}f(x){\rm d}x .
\ee

\subsection{}

First of all, we have from (\ref{3.2}), for example
\be \label{3.3}
\int_{\mg_{2\nu(-x)}}^{\mg_{2\nu(x)}}f[\vp_1(t)]\tilde{Z}^2(t){\rm d}t=\int_{g_{2\nu(-x)}}^{g_{2\nu(x)}}f(t){\rm d}t
\ee
(see (\ref{1.4})). Next, in the case $f(t)=Z^2(t)$, we have the following $\tilde{Z}^2$-transformation
\be \label{3.4}
\begin{split}
& \int_{\mG_3(x)}Z^2[\vp_1(t)]\tilde{Z}^2(t){\rm d}t=\int_{G_3(x)}Z^2(t){\rm d}t , \\
& \int_{\mG_4(y)}Z^2[\vp_1(t)]\tilde{Z}^2(t){\rm d}t=\int_{G_4(y)}Z^2(t){\rm d}t .
\end{split}
\ee
Let us remind that we have proved in the paper \cite{3} the following mean-value formulae
\be \label{3.5}
\begin{split}
& \int_{G_3(x)}Z^2(t){\rm d}t=\frac x\pi U\ln\frac{T}{2\pi}+\frac{2x}{\pi}\left( c+\frac{\sin x}{x}\right)U+\mcal{O}(xT^{5/12}\ln^2T) , \\
& \int_{G_4(y)}Z^2(t){\rm d}t=\frac y\pi U\ln\frac{T}{2\pi}+\frac{2y}{\pi}\left( c-\frac{\sin y}{y}\right)U+\mcal{O}(yT^{5/12}\ln^2T) .
\end{split}
\ee
Now, the formula (\ref{1.5}) follows from (\ref{3.4}), (\ref{3.5}).

\begin{remark}
The formulae (\ref{3.5}) are the consequences of their discrete form
\bdis
\begin{split}
& \sum_{T\leq g_{2\nu}\leq T+U}Z^2[g_{2\nu}(\tau)] = \\
& \frac{1}{2\pi}U\ln^2\frac{T}{2\pi}+\frac 1\pi(c+\cos\tau)U\ln\frac{T}{2\pi}+\mcal{O}(T^{5/12}\ln^3T) , \\
& \sum_{T\leq g_{2\nu+1}\leq T+U}Z^2[g_{2\nu+1}(\tau)] = \\
& \frac{1}{2\pi}U\ln^2\frac{T}{2\pi}+\frac 1\pi(c-\cos\tau)U\ln\frac{T}{2\pi}+\mcal{O}(T^{5/12}\ln^3T)
\end{split}
\edis
(see \cite{3}, (10)).
\end{remark}

\thanks{I would like to thank Michal Demetrian for helping me with the electronic version of this work.}

\end{document}